\documentclass[11pt]{article}
\usepackage{tikz}
\usetikzlibrary{arrows}
\tikzstyle{block}=[draw opacity=0.7,line width=1.4cm]

\usepackage{geometry}
\usepackage{graphicx}
\usepackage{setspace}
\usepackage{verbatim}
\usepackage{url}
\usepackage{latexsym}
\usepackage{amsmath,amsthm,amsfonts,amssymb,amscd,amsbsy}
\newcommand{\al}{{\alpha}}
\newcommand{\re}{{\mathbb{R}}}
\newcommand{\nb}{{\mathbb{N}}}
\newcommand{\zb}{{\mathbb{Z}}}
\newcommand{\be}{\beta}
\newcommand{\und}{\underline}

\newcommand{\R}{\mathbb R}
\newcommand{\N}{\mathbb N}
\newtheorem{theorem}{Theorem}[section]
\newtheorem{corollary}[theorem]{Corollary}
\newtheorem{lemma}[theorem]{Lemma}
\newtheorem{proposition}[theorem]{Proposition}

\newtheorem*{question*}{Question}

\theoremstyle{definition}
\newtheorem{definition}[theorem]{Definition}
\newtheorem{remark}[theorem]{Remark}

\def\cK{{\cal {K}}}
\def\sm{\smallskip} \def\ens{\enspace}
\def\df={\buildrel {\rm def}\over =}

\begin{document}

\title{Dynamical Systems, Fractal Geometry and Diophantine Approximations 
  }

\author{Carlos Gustavo Tamm de Araujo Moreira \\
IMPA}

\maketitle
\begin{abstract}
We describe in this survey several results relating Fractal Geometry, Dynamical Systems and Diophantine Approximations, including a description of recent results related to geometrical properties of the classical Markov and Lagrange spectra and generalizations in Dynamical Systems and Differential Geometry.
\end{abstract}

%\begin{classification}
%Primary 37C29, 28A80
%\end{classification}

%\begin{keywords}
%Fractal geometry. Homoclinic bifurcations. 
%\end{keywords}

%\maketitle

\section{Introduction}
The theory of Dynamical Systems is concerned with the asymptotic behaviour of systems which evolve over time, and gives models for many phenomena arising from natural sciences, as Meteorology and Celestial Mechanics. The study of a number of these models had a fundamental impact in the developement of the mathematical theory of Dynamical Systems. An important initial stage of the theory of Dynamical Systems was the study by Poincar\'e of the restricted three-body problem in Celestial Mechanics in late nineteenth century, during which he started to consider the qualitative theory of differential equations and proved results which are also basic to Ergodic Theory, as the famous Poincar\'e's recurrence lemma. He also discovered during this work the {\it homoclinic} behaviour of certain orbits, which became very important in the study of the dynamics of a system. The existence of transverse homoclinic points implies that the dynamics is quite complicated, as remarked already by Poincar\'e: ``Rien n'est plus propre \`a nous donner une id\'ee de la complication du probl\`eme des trois corps et en g\'en\'eral de tous les probl\`emes de Dynamique...", in his classic {\it Les M\'ethodes Nouvelles de la M\'ecanique C\'eleste} (\cite{Po}), written in late 19th century. This fact became clearer much decades later, as we will discuss below.

Poincar\'e's original work on the subject was awarded a famous prize in honour of the 60th birthday of King Oscar II of Sweden. There is an interesting history related to this prize. In fact, there were two versions of Poincar\'e's work presented for the prize, whose corresponding work was supposed to be published in the famous journal Acta Mathematica - a mistake was detected by the Swedish mathematician Phragm\'en in the first version. When Poincar\'e became aware of that, he rewrote the paper, including the quotation above calling the attention to the great complexity of dynamical problems related to homoclinic intersections.  We refer to the excellent paper \cite{Y} by Jean-Christophe Yoccoz describing such an event.

We will focus our discussion of Dynamical Systems on the study of flows (autonomous ordinary differential equations) and iterations of diffeomorphisms, with emphasis in the second subject (many results for flows are similar to corresponding results for diffeomorphisms). 

The first part of this work is related to the interface between Fractal Geometry and Dynamical Systems. We will give particular attention to results related to dynamical bifurcations, specially {\it homoclinic bifurcations}, perhaps the most important mechanism that creates complicated dynamical systems from simple ones. We will see how the study of the fractal geometry of hyperbolic sets has a central r\^ole in the study of dynamical bifurcations. We shall discuss recent results and ongoing works on fractal geometry of hyperbolic sets in arbitrary dimensions.

The second part is devoted to the study of the interface between Fractal Geometry and Diophantine Approximations. The main topic of this section will be the study of geometric properties of the classical Markov and Lagrange spectra - we will see how this study is related to the study of sums of regular Cantor sets, a topic which also appears naturally in the study of homoclinic bifurcations (which seems, at a first glance, to be a very distant subject from Diophantine Approximations).

The third part is related to the study of natural generalizations of the classical Markov and Lagrange spectra in Dynamical Systems and in Differential Geometry - for instance, we will discuss properties of generalized Markov and Lagrange spectra associated to geodesic flows in manifolds of negative curvature and to other hyperbolic dynamical systems. This is a subject with much recent activity, and several ongoing relevant works.

{\bf Acknowledgements:} We would like to warmly thank Jacob Palis for very valuable discussions on this work.
%\pagebreak

\section{Fractal Geometry and Dynamical Systems}

\subsection{Hyperbolic sets and Homoclinic Bifurcations}

The notion of {\it hyperbolic systems} was introduced by Smale in the sixties, after a global example provided by Anosov, namely the diffeomorphism $f(x,y)=(2x+y, x+y)\pmod 1$ of the torus ${\mathbb T}^2={\mathbb R}^2/{\mathbb Z}^2$ and a famous example, given by Smale himself, of a {\it horseshoe}, that is a robust example of a dynamical system on the plane with a transverse homoclinic point as above, which implies a rich dynamics - in particular the existence of infinitely many periodic orbits. The figure below depicts a horseshoe.

%\begin{figure}
%\centering
\begin{center}
\includegraphics[width=5cm]{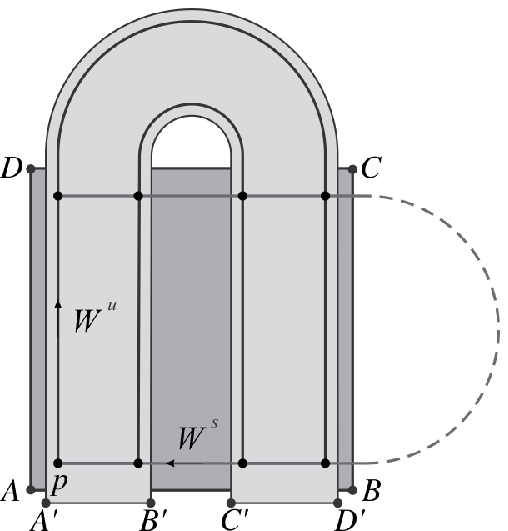}
\end{center}
%\end{figure}

(the dynamics sends the square $(ABCD)$ onto the domain bounded by $(A'B'C'D')$).

Let $\Lambda\subset M$ be a compact subset of a manifold $M$.  We say  that $\Lambda$ is a hyperbolic set for a diffeomorphism $\varphi\colon
M\to M$ if $\varphi(\Lambda) = \Lambda$ and there is a decomposition $T_\Lambda M
= E^s\oplus E^u$ of the tangent bundle of $M$ over $\Lambda$ such that $D\varphi\mid_{E^s}$ is uniformly contracting and
$D\varphi\mid_{E^u}$ is uniformly expanding. We say that $\varphi$ is hyperbolic if the limit set of its dynamics is a hyperbolic set.

It is important to notice that when a diffeomorphism has a transverse homoclinic point then its dynamics contains an non-trivial invariant hyperbolic set which is (equivalent to) a horseshoe - this explains the complicated situation discovered by Poincar\'e as mentioned above.

The importance of the notion of hyperbolicity is also related to the {\it stability conjecture} by Palis and Smale, according to which structurally stable dynamical systems are essentially the hyperbolic ones (i.e. the systems whose limit set is hyperbolic). After important contributions by Anosov, Smale, Palis, de Melo, Robbin and Robinson, this conjecture was proved in the $C^1$ topology by Ma\~n\'e (\cite{Ma}) for diffeomorphisms and, later, by Hayashi (\cite{Hay}) for flows. The stablility conjecture (namely the statement that structural stability implies hyperbolicity) is still open in the $C^k$ topology for $k\ge 2$.

As mentioned in the introduction, homoclinic bifurcations are perhaps the most important mechanism that creates complicated dynamical systems from simple ones. This phenomenon takes place when an element of a family of dynamics (diffeomorphisms or flows) presents a hyperbolic periodic point whose stable and unstable manifolds have a non-transverse intersection. When we connect, through a family, a dynamics with no {\it homoclinic points} (namely, intersections of stable and unstable manifolds of a hyperbolic periodic point) to another one with a {\it transverse} homoclinic point (a homoclinic point where the intersection between the stable and unstable manifolds is transverse) by a family of dynamics, we often go through a homoclinic bifurcation. 

Homoclinic bifurcations become important when going beyond the hyperbolic theory. In the late sixties, Sheldon Newhouse combined homoclinic bifurcations with the complexity already available in the hyperbolic theory and some new concepts in Fractal Geometry to obtain dynamical systems far more complicated than the hyperbolic ones. Ultimately this led to his famous result on the coexistence of infinitely many periodic attractors (which we will discuss later). Later on, Mora and Viana (\cite{MV}) proved that any surface diffeomorphism presenting a homoclinic tangency can be approximated by a diffeomorphism exhibiting a H\'enon-like strange attractor (and that such diffeomorphisms appear in any typical family going through a homoclinic bifurcation).

Palis conjectured that any diffeomorphism of a surface can be approximated arbitrarily well in the $C^k$ topology by a hyperbolic diffeomorphism or by a diffeomorphism displaying a homoclinic tangency. This was proved by Pujals and Sambarino (\cite{PS}) in the $C^1$ topology. Palis also proposed a general version of this conjecture: any diffeomorphism (in arbitrary ambient dimension) can be approximated arbitrarily well in the $C^k$ topology by a hyperbolic diffeomorphism, by a diffeomorphism displaying a homoclinic tangency or by a diffeomorphism displaying a heteroclinic cycle (a cycle given by intersections of stable and unstable manifolds of periodic points of different indexes). A major advance was done by Crovisier and Pujals (\cite{CP}), who proved that any diffeomorphism can be approximated arbitrarily well in the $C^1$ topology by a diffeomorphism displaying a homoclinic tangency, by a diffeomorphism displaying a heteroclinic cycle or by an {\it essentially hyperbolic} diffeomorphism: a diffeomorphism which displays a finite number of attractors whose union of basins is open and dense.

The first natural problem related to homoclinic bifurcations is the study of {\it homoclinic explosions} on surfaces: We consider one-parameter families $(\varphi_{\mu})$, $\mu \in (-1,1)$ of diffeomorphisms of a surface for which $\varphi_{\mu}$ is uniformly hyperbolic for $\mu<0$, and $\varphi_0$ presents a quadratic homoclinic tangency associated to a hyperbolic periodic point (which may belong to a {\it horseshoe} - a compact, locally maximal, hyperbolic invariant set of saddle type). It unfolds for $\mu>0$ creating locally two transverse intersections between the stable and unstable manifolds of (the continuation of) the periodic point. A main question is: what happens for (most) positive values of $\mu$? The following figure depicts such a situation for $\mu=0$.

%\begin{figure}
%\centering
\begin{center}
\includegraphics[width=5cm]{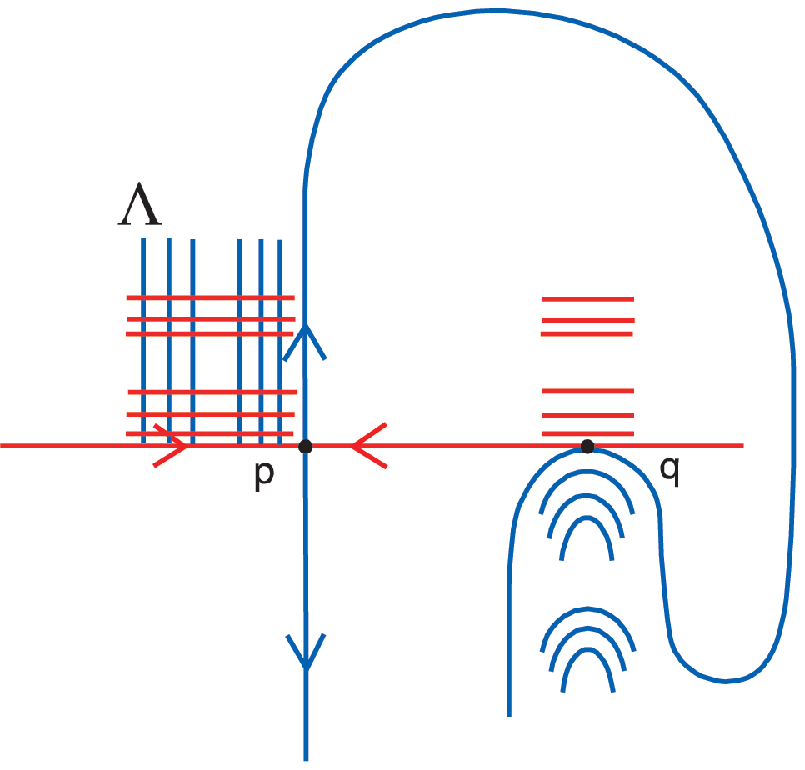}
\end{center}
%\end{figure}

Fractal sets appear naturally in Dynamical Systems and fractal dimensions when we try to measure fractals. They are essential to describe most of the main results in this presentation.

Given a metric space $X$, it is often true that the minimum number $N(r)$ of balls of radius $r$ needed to cover $X$ is roughly proportional to $1/r^d$, for some positive constant $d$, when $r$ becomes small. In this case, $d$ will be the {\it box dimension} of $X$. More precisely, we define the (upper) box dimension of $X$ as 
$$d(X)=\limsup_{r\to 0} \frac{\log N(r)}{-\log r}.$$

The notion of {\it Hausdorff dimension} of a set is more subtle, but more useful. The main difference with the notion of box dimension is that, while the box dimension is related to coverings of $X$ by small balls of equal radius, the Hausdorff dimension deals with arbitrary coverings of $X$ by balls of small (but not necessarily equal) radius. 

Given a countable covering $\cal U$ of $X$ by
balls, ${\cal U} = (B(x_i,r_i))_{i\in \N}$, we define its norm $||{\cal U}||$ as 
$||{\cal U}|| =
\max\{r_i,\,i\in \N\}$ (where $r_i$ is the radius of the ball $B(x_i,r_i)$). Given
$s\in\re_+$, we define ${\cal H}_s({\cal U}) =\sum_{i\in \N} r_i^s$.

The Hausdorff $s$-measure of $X$ is 
$${\cal H}_s(X) = \lim_{\epsilon\to
0} \inf_{{\cal U} \text{ covers } X \atop ||{\cal U}||<\epsilon}{\cal H}_s({\cal U}).$$

One can show that there is an unique real number, the {\it Hausdorff
dimension} of $X$, which we denote by $HD(X)$, such that $s < HD(X) \Rightarrow
{\cal H}_s(X) = +\infty$ and $s > HD(X) \Rightarrow {\cal H}_s(X) = 0$ (so $HD(X)$ can be defined shortly as
$$HD(X)=\inf \{s>0; \inf_{X\subset \cup B(x_n,r_n)}\sum r_n^s=0\}).$$

For ``well-behaved" sets $X$ - in particular for regular Cantor sets in ambient dimension $1$ and horseshoes in ambient dimension $2$, the box and Hausdorff dimensions of $X$ coincide.
    
Regular Cantor sets on the line play a fundamental role in dynamical systems and notably also in some problems in number theory. They are defined by expanding maps and have some kind of self-similarity property: small parts of them are diffeomorphic to big parts with uniformly bounded distortion (we will give a precise definition in a while). In both settings, dynamics and number theory, a key question is whether the arithmetic difference (see definition below) of two such sets has non-empty interior.

A horseshoe $\Lambda$ in a surface is locally diffeomorphic to the Cartesian product of two regular Cantor sets: the so-called {\it stable} and {\it unstable} Cantor sets $K^s$ and $K^u$ of $\Lambda$, given by intersections of $\Lambda$ with local stable and unstable manifolds of some points of the horseshoe. The Hausdorff dimension of $\Lambda$, which is equal to the sum of the Hausdorff dimensions of $K^s$ and $K^u$, plays a fundamental role in several results on homoclinic bifurcations associated to $\Lambda$.

From the dynamics side, in the eighties, Palis and Takens (\cite{PT}, \cite{PT1})
proved the following theorem about homoclinic bifurcations associated to a hyperbolic set:

\begin{theorem}
Let $(\varphi_{\mu})$, $\mu \in (-1,1)$ be a family of diffeomorphisms of a surface presenting a homoclinic explosion at $\mu=0$ associated to a periodic point belonging to a horseshoe $\Lambda$. Assume that $HD(\Lambda)<1$. Then
$$\lim_{\delta \to 0}\frac{m(H\cap [0,\delta])}{\delta}=1,$$
where $H:=\{\mu>0 \mid \varphi_{\mu}\,\hbox{  is uniformly hyperbolic}\}$.
\end{theorem}

\subsection{Regular Cantor sets - a conjecture by Palis}

A central fact used in the proof of the above theorem by Palis and Takens is that if $K_1$ and $K_2$ are regular Cantor sets on the real line such that the sum of their Hausdorff dimensions is smaller than one, then $K_1-K_2 =
\{x-y \mid x \in K_1, y \in K_2\}=\{t\in \R|K_1\cap(K_2+t)\ne \emptyset\}$ (the {\it arithmetic difference} between $K_1$ and $K_2$) is a set of zero Lebesgue measure (indeed of Hausdorff dimension smaller than 1). On the occasion, looking for some kind of characterization property for this phenomenon, Palis conjectured (see \cite{P}, \cite{P1}) that for generic pairs of regular Cantor sets $(K_1, K_2)$ of the real line either $K_1- K_2$ has zero measure or else it contains an interval (the last case should correspond in homoclinic bifurcations to open sets of tangencies). A slightly stronger statement is that,
 if $K_1$ and $K_2$ are generic regular Cantor sets and the sum of their Hausdorff dimensions is bigger than 1, then $K_1-K_2$ contains intervals. 

Another motivation for the conjecture was Newhouse's work in the seventies, when he introduced the concept of thickness of a regular Cantor set, another fractal invariant associated to Cantor sets on the real line. It was used in \cite{N1} to exhibit open sets of
diffeomorphisms with persistent homoclinic tangencies, therefore with no
hyperbolicity. It is possible (\cite{N2}) to prove that, under a dissipation hypothesis, in such an open set there is
a residual set of diffeomorphisms which present infinitely many coexisting
sinks. In \cite{N3}, it is proved that under generic hypotheses every family of
surface diffeomorphisms that unfold a homoclinic tangency goes through such an
open set.  It is to be noted that in the case described above with $HD(\Lambda)<1$ (as studied in \cite{PT1}) these sets have zero density.  See \cite{PT2} for a detailed presentation of these results. An important related question by Palis is whether the sets of parameter values corresponding to infinitely many coexisting sinks have typically zero Lebesgue measure.

An earlier and totally independent development had taken place in number theory. In 1947, M. Hall (\cite{H}) proved that any real number can be written as the sum of two numbers whose continued fractions  coefficients (of positive index) are bounded by $4$. More precisely, if $C(4)$ is the regular Cantor set (see general definition below) formed of such numbers in $[0,1]$, then one has $C(4)+C(4) = [\sqrt2 -1, 4(\sqrt 2 -1)]$. We will discuss generalizations and consequences of this result in the next section.

A Cantor set $K\subset {\mathbb R}$ is a {\it $C^k$-regular Cantor set}, $k\ge 1$, if:

\begin{itemize}

\item[i)] there are disjoint compact intervals $I_1,I_2,\dots,I_r$ such that $K\subset I_1 \cup \cdots\cup I_r$ and the boundary of each $I_j$ is contained in $K$;

\item[ii)] there is a $C^k$ expanding map $\psi$ defined in a neighbourhood of $I_1\cup I_2\cup\cdots\cup I_r$ such that, for each $j$, $\psi(I_j)$ is the convex hull of a finite union of some of these intervals $I_s$. Moreover, we suppose that $\psi$ satisfies:

\begin{itemize}

\item[ii.1)] for each $j$, $1\le j\le r$ and $n$ sufficiently big, $\psi^n(K\cap I_j)=K$;

\item[ii.2)] $K=\bigcap\limits_{n\in\mathbb N} \psi^{-n}(I_1\cup I_2\cup\cdots\cup I_r)$.

\end{itemize}
\end{itemize}

\begin{remark}
If $k$ is not an integer, say $k=m+\alpha$, with $m\ge 1$ integer and $0<\alpha<1$ we assume that $\psi$ is $C^m$ and $\psi^{(m)}$ is $\alpha$-H\"older.
\end{remark}

We say that $\{I_1,I_2,\dots,I_r\}$ is a {\it Markov partition} for $K$ and that $K$ is {\it defined} by $\psi$.

\begin{remark}
In general, we say that a set $X \subset \R$ is a Cantor set if $X$ is compact, without isolated points and with empty interior. Cantor sets in $\R$ are homeomorphic to the classical ternary Cantor set $K_{1/3}$ of the elements of $[0,1]$ which can be written in base 3 using only digits $0$ and $2$. The set $K_{1/3}$ is itself a regular Cantor set, defined by the map $\psi:[0,1/3]\cup [2/3,1] \to \R$ given by $\psi[x]=3x-\lfloor 3x \rfloor$.
\end{remark}

An {\it interval of the construction of the regular Cantor set $K$} is a connected component of $\psi^{-n}(I_j)$ for some $n\in\mathbb N$, $j\le r$.

Given $s \in [1,k]$ and another regular Cantor set $\tilde K$, we say that $\tilde K$ is close to $K$ in the $C^s$ topology if $\tilde K$ has a Markov partition $\{\tilde I_1,\tilde I_2,\dots,\tilde I_r\}$ such that the interval $\tilde I_j$ has endpoints close to the endpoints of $I_j$, for $1 \le j \le r$ and $\tilde K$ is defined by a $C^s$ map $\tilde \psi$ which is close to $\psi$ in the $C^s$ topology.

The $C^{1+}$-topology is such that a sequence $\psi_n$ converges to $\psi$ if there is some $\alpha>0$ such that $\psi_n$ is $C^{1+\alpha}$ for every $n\ge 1$ and $\psi_n$ converges to $\psi$ in the $C^{1+\alpha}$-topology.

The concept of stable intersection of two regular
Cantor sets was introduced in \cite{M}: two Cantor sets $K_1$ and $K_2$ have stable intersection if there is a neighbourhood $V$ of $(K_1,K_2)$ in the set of pairs of $C^{1+}$-regular Cantor sets such that $(\widetilde
K_1, \widetilde K_2) \in V \Rightarrow \widetilde K_1 \cap \widetilde K_2 \ne \emptyset$. 

In the same paper conditions based on
renormalizations were introduced to ensure stable intersections, and applications of stable intersections to homoclinic bifurcations were obtained: roughly speaking, if some translations of the stable and unstable regular Cantor sets associated to the horseshoe at the initial bifurcation parameter $\mu=0$ have stable intersection then the set $\{\mu>0 \mid \varphi_{\mu}\,\,\hbox{ presents persistent homoclinic tangencies}\}$ has positive Lebesgue density at $\mu=0$. It was also shown that this last phenomenon can coexist with positive density of hyperbolicity in a
persistent way. 

Besides, the following question was posed in \cite{M}:
Does there exist a dense (and automatically open) subset ${\cal U}$ of
$$\Omega^{\infty} =\{(K_1,K_2); K_1, K_2\,\,\hbox{are} \,C^{\infty}-\hbox{regular Cantor sets and}\, HD(K_1) + HD(K_2) > 1\}$$
%$$\Omega^{\infty} =\{(K_1,K_2); K_1, K_2\,\,\,\, C^{\infty}-\hbox{regular Cantor sets}\mid HD(K_1) + HD(K_2) > 1\}$$
such that $(K_1,K_2) \in {\cal U} \Rightarrow \exists\, t\in\re$ such that $(K_1,K_2+t)$ has stable intersection?
A positive answer to this question implies a strong version of Palis'
conjecture. Indeed, $K_1-K_2=\{t\in \R\mid K_1\cap (K_2+t)\ne \emptyset\}$, so, if $(K_1,K_2+t)$ has stable intersection then $t$ belongs persistently to the interior of $K_1-K_2$.

Moreira and Yoccoz \cite{MY} gave an affirmative answer to this question, proving the following

\begin{theorem}
There is an open and dense set \,\,${\cal U} \subset \Omega^{\infty}$ such that if $(K_1,K_2) \in \cal{U}$, then $I_s(K_1,K_2)$ is dense in $K_1-K_2$ and $HD((K_1-K_2)\backslash I_s(K_1,K_2)) < 1$, where 
$$I_s(K_1,K_2) := \{t \in \mathbb{R} \mid K_1\hbox{ and }(K_2+t)\hbox{ have stable intersection}\}.$$
\end{theorem}

\begin{center}
\includegraphics[width=10cm]{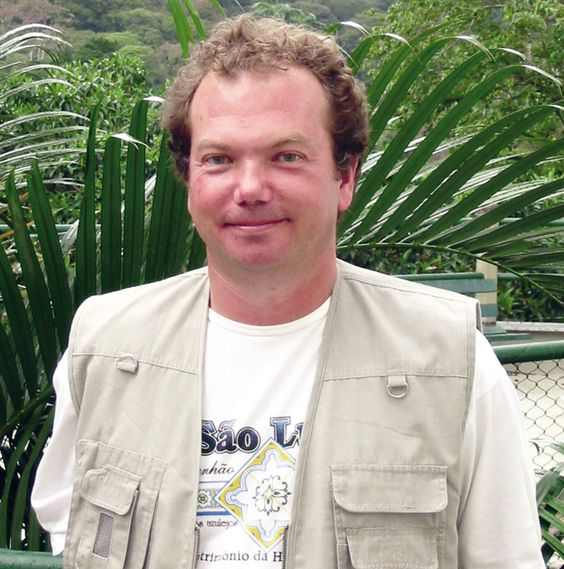}
\end{center}
\hskip 5.5cm {\bf Jean-Christophe Yoccoz}
\vskip .2in

The same result works if we replace stable intersection
by $d$-stable intersection, which is defined by asking that any pair $(\widetilde K_1, \widetilde K_2)$ in some neighbourhood of $(K_1,K_2)$  satisfies $HD(\widetilde K_1 \cap \widetilde K_2) \ge d$: most pairs of Cantor sets $(K_1,K_2) \in {\Omega}^{\infty}$ have $d$-stable intersection for any $d < HD(K_1) + HD(K_2)-1$.

The open set $\cal U$ mentioned in the above theorem is very large in ${\Omega}^\infty$ in the sense that generic $n$-parameter families in ${\Omega}^\infty$ are actually contained in $\cal U$.

The proof of this theorem depends on a sufficient condition for the existence of stable intersections of two Cantor sets, related to a notion of renormalization, based on the fact that small parts of regular Cantor sets are diffeomorphic to the whole set: the existence of a {\it recurrent compact set} of relative positions of limit geometries of them. Roughly speaking, it is a compact set of relative positions of regular Cantor sets such that, for any relative position in such a set, there is a pair of (small) intervals of the construction of the Cantor sets such that the renormalizations of the Cantor sets associated to these intervals belong to the interior of the same compact set of relative positions. 

The main result is reduced to prove the existence of recurrent compact sets of relative positions for most pairs of regular Cantor sets whose sum of Hausdorff dimensions is larger than one. A central argument in the proof of this fact is a probabilistic argument \`a la Erd\H os: we construct a family of perturbations with a large number of parameters and show the existence of such a compact recurrent set with large probability in the parameter space (without exhibiting a specific perturbation which works). 

Another important ingredient in the proof is the Scale Recurrence Lemma, which, under mild conditions on the Cantor sets (namely that at least one of them is not essentially affine), there is a recurrent compact set for renormalizations at the level of relative scales of limit geometries of the Cantor sets. This lemma is the fundamental tool in the paper \cite{M4}, in which it is proved that, under the same hypothesis above, if $K$ and $K'$ are $C^2$-regular Cantor sets, then $HD(K+K')=\min\{1, HD(K)+HD(K')\}$.

An important result in fractal geometry which is used in \cite{MY} is the famous Marstrand's theorem (\cite{Mars}), according to which, given a Borel set $X \subset {\mathbb R}^2$ with $HD(X)>1$ then, for almost every $\lambda \in {\mathbb R}$, $\pi_{\lambda}(X)$ has positive Lebesgue measure, where $\pi_{\lambda}:{\mathbb R}^2 \to {\mathbb R}$ is given by $\pi_{\lambda}(x,y)=x-\lambda y$. In particular, if $K_1$ and $K_2$ are regular Cantor sets with $HD(K_1)+HD(K_2)>1$ then, for almost every $\lambda \in {\mathbb R}$, $K_1-\lambda K_2$ has positive Lebesgue measure. Moreira and Yuri Lima gave combinatorial alternative proofs of Marstrand's theorem, first in the case of Cartesian products regular Cantor sets (\cite{LiM1}) and then in the general case (\cite{LiM2}).

In \cite{MY2}, Moreira and Yoccoz proved the following fact concerning
generic homoclinic bifurcations associated to two dimensional 
saddle-type hyperbolic sets (horseshoes) with Hausdorff dimension bigger than one: typically there are translations of the stable and unstable Cantor sets having stable intersection, and so it yields open sets of stable tangencies in the parameter line with positive density at the initial bifurcation value. Moreover, the union of such a set with the hyperbolicity set in the parameter line generically has full density at the initial bifurcation value. This extends the results of \cite{PY4}.

The situation is quite different in the $C^1$-topology, in which stable intersections do not exist:

\begin{theorem}[\cite{M2}]
Given any pair $(K,K')$ of  regular Cantor sets, we can find, arbitrarily close to it in the $C^1$ topology, pairs $(\tilde K,\tilde K')$ of regular Cantor sets with $\tilde K\cap \tilde K'=\emptyset$.

Moreover, for generic pairs $(K,K')$ of $C^1$-regular Cantor sets, the arithmetic difference $K-K'$ has empty interior (and so is a Cantor set).
\end{theorem}

The main technical difference between the $C^1$ case and the $C^2$ (or even $C^{1+\al}$) cases is the lack of bounded distortion of the iterates of $\psi$ in the $C^1$ case, and this fact is fundamental for the proof of the previous result.

The previous result may be used to show that there are no $C^1$ robust tangencies between leaves of the stable and unstable foliations of respectively two given hyperbolic horseshoes $\Lambda_1, \Lambda_2$ of a diffeomorphism of a surface. This is also very different from the situation in the $C^\infty$ topology - for instance, in \cite{MY2} it is proved that, in the unfolding of a homoclinic or heteroclinic tangency associated to two horseshoes, when the sum of the correspondent stable and unstable Hausdorff dimensions is larger than one, there are generically stable tangencies associated to these two horseshoes. This result is done in the following

\begin{theorem}[\cite{M2}]
Given a $C^1$ diffeomorphism $\psi$ of a surface $M$ having two (non necessarily disjoint) horseshoes $\Lambda_1, \Lambda_2$, we can find, arbitrarily close to it in the $C^1$ topology, a diffeomorphism $\tilde \psi$ of the surface for which the horseshoes $\Lambda_1, \Lambda_2$ have hyperbolic continuations $\tilde \Lambda_1, \tilde \Lambda_2$, and there are no tangencies between leaves of the stable and unstable foliations of $\tilde \Lambda_1$ and $\tilde \Lambda_2$, respectively.
Moreover, there is a generic set $\cal R$ of $C^1$ diffeomorphisms of $M$ such that, for every $\check \psi \in \cal R$, there are no tangencies between leaves of the stable and unstable foliations of $\Lambda_1, \Lambda_2$, for any horseshoes $\Lambda_1, \Lambda_2$ of $\check \psi$.
\end{theorem}

Since stable intersections of Cantor sets are the main known obstructions to density of hyperbolicity for diffeomorphisms of surfaces, the previous result gives some hope of proving density of hyperbolicity in the $C^1$ topology for diffeomorphisms of surfaces, a well-known question by Smale. In particular in the work \cite{MMP} on a family of two-dimensional maps (the so-called Benedicks-Carleson {\it toy model} for H\'enon dynamics) in which we prove that in this family there are diffeomorphisms which present stable homoclinic tangencies (Newhouse's phenomenon) in the $C^2$-topology but their elements can be arbitrarily well approximated in the $C^1$-topology by hyperbolic maps. 

\subsection{Palis-Yoccoz theorem on non-uniformly hyperbolic horseshoes}

In a major work by Palis and Yoccoz (\cite{PY8}; see also \cite{PY6} and \cite{PY7}), they propose to advance considerably the current knowledge on the topic of bifurcations of heteroclinic cycles for smooth ($C^{\infty}$) parametrized families $\{g_t, t\in \R\}$ of surface diffeomorphisms. They assume that $g_t$ is hyperbolic for $t<0$ and $|t|$ small and that a quadratic tangency $q$ is formed at $t=0$ between the stable and unstable lines of two periodic points, not belonging to the same orbit, of a (uniformly hyperbolic) horseshoe $\Lambda$ and that such lines cross each other with positive relative speed as the parameter evolves, starting at $t=0$ near the point $q$. They also assume that, in some neighbourhood $W$ of the union of $\Lambda$ with the orbit of tangency $O(q)$, the maximal invariant set for $g_0=g_{t=0}$ is $\Lambda \cup O(q)$, where $O(q)$ denotes the orbit of $q$ for $g_0$.

%\begin{figure}
%\centering
\begin{center}
\includegraphics[width=12cm]{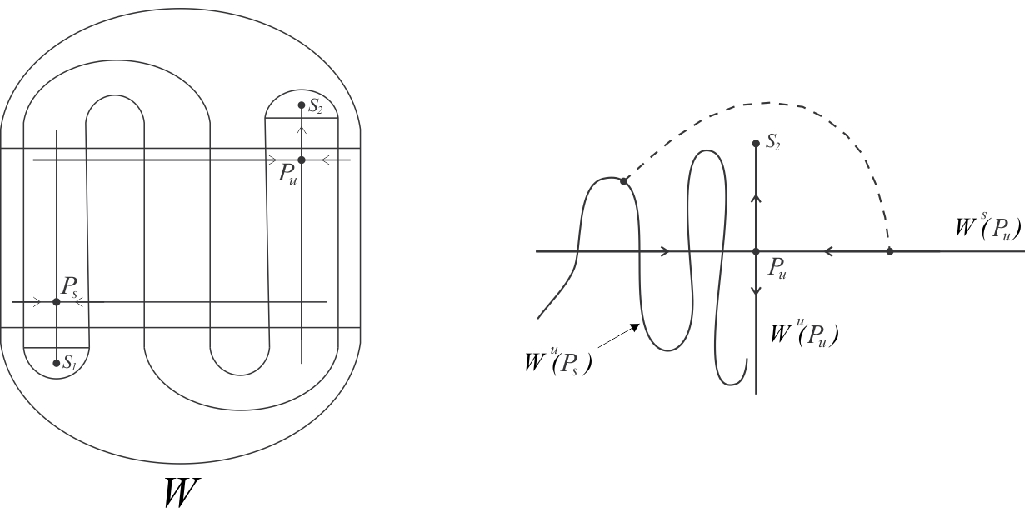}
\end{center}
%\end{figure}

A main novelty is that they allowed the Hausdorff dimension $HD(\Lambda)$ to be bigger than one: they assume $HD(\Lambda)$ to be bigger than one, but not much bigger (more precisely, they assume that, if $d_s$ and $d_u$ are the Hausdorff dimensions of, respectively, the stable and unstable Cantor sets of $g_0$ then $(d_s+d_u)^2+\max(d_s, d_u)^2<d_s+d_u+\max(d_s, d_u)$). Then, for most small values of $t$, $g_t$ is a ``non-uniformly hyperbolic" horseshoe in $W$, and so $g_t$ has no attractors nor repellors in $W$. Most small values of $t$, and thus most $g_t$, here means that $t$ is taken in a set of parameter values with Lebesgue density one at $t=0$.   

The construction of non-uniformly hyperbolic horseshoes for most parameters is a highly non-trivial counterpart of Yoccoz' proof \cite{Y-J} (based on the so-called \emph{Yoccoz puzzles}) of the celebrated Jakobson's theorem in the context of heteroclinic explosions. 

Of course, Palis and Yoccoz do not consider their result as the end of the line. Indeed, they expected the same results to be true for all cases $0<HD(\Lambda)<2$. However, to achieve that, it seems that their methods need to be considerably sharpened: it would be necessary to study deeper the dynamical recurrence of points near tangencies of higher order (cubic, quartic...) between stable and unstable curves. They also expected their results to be true in higher dimensions.

Finally, they hoped that the ideas introduced in that work might be useful in broader contexts. In the horizon lies a famous question concerning the standard family of area preserving maps (which is the family $(f_k)_{k \in \R}$ of diffeomorphisms of the torus ${\mathbb T}^2=\R^2/{\mathbb Z}^2$ given by $f_k(x, y)=(-y+2x+k\sin(2\pi x), x) \pmod 1$): can we find sets of positive  Lebesgue measure in the parameter space such that the corresponding maps display non-zero Lyapunov exponents in sets of positive Lebesgue probability in phase space?

In a recent work, Matheus, Moreira and Palis applied the results of \cite{PY8} to prove the following

\begin{theorem}\label{t.MMP} There exists $k_0 > 0$ such that, for all $|k|>k_0$ and $\varepsilon>0$, the subset of parameters $r\in\mathbb{R}$ such that $|r-k|<\varepsilon$ and $f_r$ exhibits a non-uniformly hyperbolic horseshoe (in the sense of Palis--Yoccoz \cite{PY8}) has positive Lebesgue measure. 
\end{theorem}

\subsection{Hyperbolic sets in higher dimensions}

	The fractal geometric properties of horseshoes of surface diffeomorphisms are quite regular and well-understood. As we mentioned before, a horseshoe $\Lambda$ in ambient dimension $2$ is locally diffeomorphic to a cartesian product $K^s \times K^u$, where $K^s$ and $K^u$ are regular Cantor sets on the real line. McCluskey and Manning proved that $HD(\Lambda)$ is continuous in the $C^1-$topology (later Palis and Viana gave a more geometrical proof of this fact). 

 In ambient dimension larger than two, much less is known about the geometry of horseshoes. In general, if $\Lambda$ is such a horseshoe, the Hausdorff dimension $HD(\Lambda)$ is not always continuous, as shown by Bonatti, D\'iaz and Viana. A natural question is whether $HD(\Lambda)$ is ``typically" continuous. Another natural question is what can be said about the fractal geometry of stable and unstable Cantor sets of horseshoes in arbitrary dimensions. These Cantor sets are no more conformal sets with bounded geometry, which was the case of stable and unstable Cantor sets of horseshoes in surfaces. In what follows, we will discuss some recent results on the geometry of horseshoes in arbitrary dimensions.
	
	 In the first work, in collaboration with J. Palis and M. Viana, we generalize an important lemma from \cite{MPV}: given a horseshoe $\Lambda$ whose stable spaces have dimension $\ell$ we define a family of fractal dimensions (the so-called {\it upper stable dimensions}) $\overline d^{(j)}_s(\Lambda), 1\le j\le \ell$ which satisfy $\overline d^{(1)}_s(\Lambda) \ge \overline d^{(2)}_s(\Lambda) \ge \dots \ge \overline d^{(\ell)}_s(\Lambda) \ge HD(\Lambda \cap W^s(x)), \forall x \in \Lambda$. We define $\hat{d}_s (\Lambda)=\bar{d}^{(\ell)}_s (\Lambda)$. When $r<\overline d^{(r)}_s(\Lambda)\le r+1$ we have $\hat{d}_s (\Lambda)=\overline d^{(r)}_s(\Lambda)$
These definitions are inspired in the {\it affinity dimensions}, introduced by Falconer.  
We have analogous definitions for upper unstable dimensions. We prove the following results about these dimensions: given $1\le r\le \ell$ and $\varepsilon>0$ there is a $\varepsilon-$small $C^{\infty}$ perturbation of the original diffeomorphism for which the hyperbolic continuation of $\Lambda$ has a subhorseshoe $\tilde \Lambda$ which has strong-stable foliations of codimensions $j$ for $1\le j\le r$ and which satisfies $\overline d^{(r)}_s(\tilde \Lambda)>\overline d^{(r)}_s(\Lambda)-\varepsilon$.

In the second work in progress, in collaboration with W. Silva (which extends a previous joint work in codimension $1$ - see \cite{MS}), we prove that if a horseshoe $\Lambda$ has strong stable foliations of codimensions $j$ for $1\le j\le r$ and satisfies $\hat{d}_s (\Lambda)>r$ (which is equivalent to $\overline d^{(r)}_s(\Lambda)>r$) then it has a small $C^{\infty}$ perturbation which contains a blender of codimension $r$: in particular $C^1$ images of stable Cantor sets of it (of the type $\Lambda \cap W^s(x)$) in ${\mathbb R}^r$ will typically have persistently non-empty interior. We also expect to prove that the Hausdorff dimension of these stable Cantor sets typically coincide with $\hat{d}_s (\Lambda)$, and this dimension depends continuously on $\Lambda$ on these assumptions, which would imply typical continuity of Hausdorff dimensions of stable and unstable Cantor sets of horseshoes. 

	Let us recall two main ingredients in the proof by Moreira and Yoccoz of Palis' conjecture:
\begin{itemize}
\item A {\it recurrent compact set criterion} for stable intersections (which implies that arithmetic differences persistently contain intervals).

\item An application of Erd\H os probabilistic method: a family of $C^{\infty}$ small perturbations of a regular Cantor set (the second Cantor set is fixed) with a large number of parameters such that for most parameters there is a recurrent compact set for the corresponding pair of Cantor sets.
\end{itemize}

A variation of these ingredients is present in the proof of these results in collaboration of W. Silva: 
we develop a recurrent compact set criterion which implies that a given horseshoe (with a strong-stable foliation of codimension $r$) is a $r$-codimensional blender, and we prove that, if a horseshoe has a  strong-stable foliation of codimension $r$ and satisfies $\hat{d}_s (\Lambda)>r$ (which is equivalent to $\overline d^{(r)}_s(\Lambda)>r$), then there is a small $C^{\infty}$ perturbation of it which has a recurrent compact set. In order to do this, we use the probabilistic method: we construct a family of perturbations with a large number of parameters, and show that for most parameters there is such a recurrent compact set.

An important geometrical tool in the proof of these results is the following generalization of Marstrand's theorem, which was proved in works by L\'opez, Moreira, Roma\~na and Silva:

\begin{theorem}
Let $X$ be a compact metric space, $(\Lambda, \mathcal P)$ a probability space and $\pi:\Lambda\times X\to {\mathbb R}^r$ a measurable function. Informally, one can think of $\pi_\lambda(\cdot) = \pi(\lambda,\cdot)$ as a family of projections parameterized by $\lambda$. We assume that for some positives real numbers $\alpha$ and $C$ the following transversality property is satisfied:
\begin{equation}\label{transversality}
\mathcal P[\lambda\in \Lambda: d(\pi_{\lambda}(x_1),\pi_{\lambda}(x_2))\leq\delta d(x_1,x_2)^{\alpha}]\leq C\delta^r
\end{equation}
for all $\delta>0$ and all $x_1,x_2\in X$. Assume that $\dim X>\alpha r$. Then $Leb(\pi_\lambda(X))>0$ for a.e. $\lambda\in\Lambda$ and $\int_{\Lambda}Leb(\pi_\lambda(X))^{-1}d\mathcal P<+\infty$.
\end{theorem}

\section{Fractal Geometry and Diophantine Approximations}
\subsection{The classical Markov and Lagrange spectra}
The results discussed in the previous section on regular Cantor sets have, somewhat surprisingly, deep consequences in number theory, which we discuss below. We begin by introducing the classical Lagrange spectrum.

Let $\al$ be an irrational number. According to Dirichlet's theorem, the
inequality \linebreak $|\al-\frac pq|<\frac1{q^2}$ has infinitely many
rational solutions $\frac pq$. Hurwitz improved this result by proving that
$|\al-\frac pq|<\frac1{\sqrt 5 q^2}$ also has infinitely many rational
solutions $\frac pq$ for any irrational $\al$, and that $\sqrt 5$ is the largest
constant that works for any irrational $\al$. However, for particular values of
$\al$ we can improve this constant.

More precisely, we define $k(\al):=\sup\{k>0\mid  |\al-\frac pq|<\frac 1{
kq^2}$ has infinitely many rational solutions $\frac pq\}=
\limsup_{p,q\to+\infty}\, (q|q \al-p|)^{-1}$. We have
$k(\al)\ge\sqrt 5$, $\forall \al\in\re\setminus {\mathbb Q}$ and $k\left(
\frac{1+\sqrt 5}2 \right)=\sqrt 5$. We will consider the set $L=\{k(\al) \mid
\al\in\re\setminus {\mathbb Q}$, $k(\al)<+\infty\}$.

This set is called the Lagrange spectrum. Hurwitz's theorem determines the smallest
element of $L$, which is $\sqrt 5$. This set $L$ encodes many diophantine properties of real
numbers. It is a classical subject the study of the geometric structure of $L$.
Markov proved in 1879 (\cite{Mark}) that
$$
L\cap(-\infty, 3)=\{k_1=\sqrt 5< k_2=2 \sqrt 2<k_3= \frac{\sqrt{221}}5<\dots\}
$$
where $k_n$ is a sequence (of irrational numbers whose squares are rational) converging to $3$.

The elements of the Lagrange spectrum which are smaller than $3$ are exactly the numbers of the form $\sqrt{9-\frac4{z^2}}$ where $z$ is a positive integer for which there are other positive integers $x, y$ such that $1\le x\le y\le z$ and $(x,y,z)$ is a solution of the  
{\it Markov equation}
$$x^2+y^2+z^2=3xyz.$$
Since this is a quadratic equation in $x$ (resp. in $y$), whose sum of roots is $3yz$ (resp. $3xz$), given a solution $(x,y,z)$ with $x\le y\le z$, we also have the two other solutiuons $(y,z,3yz-x)$ and $(x,z,3xz-y)$ - it is possible to prove that all solutions of the Markov equation appear in the following tree:
\begin{center}
\begin{tikzpicture}[scale=0.7, line width=1.5pt]
\tikzstyle{level 3}=[sibling distance=6cm]
\tikzstyle{level 4}=[sibling distance=2.9cm]
%\tikzstyle{level 5}=[sibling distance=2cm]
\node  {$(1,1,1)$}
 child[->] {node {$(1,1,2)$}
 child[->]  {node{$(1,2,5)$}
   child[->] {node {$(1,5,13)$}
     child[->] {node {$(1,13,34)$}
%       child[->] {node{$(1,34,89)$}}
%       child[->] {node{$(13,34,1385)$}}
}
     child[->] {node{$(5,13,194)$}}
}
   child[->] {node{$(2,5,29)$}
      child[->] {node{$(2,29,169)$}}
      child[->] {node{$(5,29,433)$}}
}
}
};
\end{tikzpicture}
\end{center}

An important open problem related to Markov's equation is the
 {\it Unicity Problem}, formulated by Frobenius about 100 years ago: for any positive integers $x_1,x_2,y_1,y_2,z$ with $x_1 \le y_1 \le z$
and $x_2 \le y_2 \le z$ such that $(x_1,y_1,z)$ and $(x_2,y_2,z)$ are solutions of Markov's equation we always have $(x_1,y_1)=(x_2,y_2)$?

If the Unicity Problem has an affirmative answer then, for every real $t<3$, its pre-image $k^{-1}(t)$ by the function $k$ above consists of a single $GL_2({\mathbb Z})$-equivalence class
(this equivalence relation is such that 
$$\alpha\sim \frac{a\alpha+b}{c\alpha+d}, \forall a, b, c, d \in {\mathbb Z}, |ad-bc|=1.)$$

Despite the ``beginning'' of the set $L$ is discrete, this is not true for the whole set $L$. As we mentioned in the introduction, M.
Hall proved in 1947 (\cite{H}) that if $C(4)$ is the regular Cantor set formed by the numbers in $[0,1]$ whose coefficients in the continued fractions expansion are bounded by $4$, then one has $C(4)+C(4) = [\sqrt2 -1, 4(\sqrt 2 -1)]$. This implies that $L$ contains a whole half line (for instance $[6,+\infty))$,
and G. Freiman determined in 1975 (\cite{F3}) the biggest half line that is contained in
$L$, which is $[c_F, +\infty)$, with
$$
c_F=\frac{2221564096 + 283748\sqrt{462}}{491993569} \cong 4,52782956616\dots\,.
$$
These last two results are based on the study of sums of regular Cantor
sets, whose relationship with the Lagrange spectrum will be explained below.

Sets of real numbers whose continued fraction representation has bounded coefficients with some combinatorial constraints, as $C(4)$, are often regular Cantor sets, which we call Gauss-Cantor sets (since they are defined by restrictions of the Gauss map $g(x)=\{1/x\}$ from $(0,1)$ to $[0,1)$ to some convenient union of intervals).

We represent below the graphics of the Gauss map $g(x)=\{\frac 1x\}$.\nopagebreak[4]
\vspace{3mm}\\
{
\small
% grafico de y=g(x)
\begin{tikzpicture}[scale=4.0, line width=1.5pt]
\draw[->] (-0.1,0)--(1.1,0);
\draw[->](0,-0.1)--(0,1.1);
\foreach \k in{1,2,...,10}{
\draw[smooth, samples=30, line width=0.3pt,domain=\k:\k+1] plot({1/\x},{\x-\k});
\draw[fill=white, line width=0.2pt] (1/\k,1) circle (0.3pt);
}
\foreach \k in{11,13,...,121}{
\draw[line width=0.3pt] (1/\k,0)--(1/\k-1/\k^2,1);
\draw[fill=white, line width=0.2pt] (1/\k,1) circle (0.3pt);
}
\fill[white](1,1) circle(0.5pt);
\node[right] at (0.6,0.666) {$y=g(x)=\left\{\frac 1x\right\}$};
\end{tikzpicture}
}
\vskip .15in
If the continued fraction of $\alpha$ is 
$$
\alpha = [a_0;a_1,a_2,\dots] \df = a_0 +\cfrac1{a_1 + \cfrac1{a_2 + {\atop\ddots}}}.
$$
 then we have the
following formula for $k(\alpha)$:
$$
k(\alpha)=\limsup_{n\to\infty}(\alpha_n+\beta_n),
$$
$$
\text{where } \alpha_n=[a_n;a_{n+1},a_{n+2},\dots]\text{ and }
\beta_n=[0;a_{n-1}, a_{n-2},\dots,a_1].
$$
The previous formula follows from the equality 
$$
 |\alpha-\frac{p_n}{q_n}| =\frac1{(\alpha_{n+1}+\beta_{n+1}) q_n^2},\quad \forall
n\in {\mathbb N},
$$
where 
$$p_n/q_n=[a_0;a_1,a_2,\dots,a_n]=a_0 +\cfrac1{a_1 + \cfrac1{a_2 + {\atop\ddots+\frac1{a_n}}}}, n \in \mathbb N$$
 are the convergents of the continued fraction of $\alpha$.

There are many results which relate the dynamics of the Gauss map with the behaviour of continued fractions. For instance, the Khintchine-L\'evy theorem, which follows from techniques of Ergodic Theory, states that, for (Lebesgue) almost every $\alpha\in\re$,
$$\lim_{n\to\infty}\sqrt[n]{q_n}=e^{\pi^2/12\log 2}.$$

{\bf Remark:}
The following elementary general facts on Diophantine approximations of real numbers show that the best rational approximations of a given real number are given by convergents of its continued fraction representation:

\noindent $\bullet$ For every $n\in\mathbb N$, 
$$\left| \alpha - \frac{p_n}{q_n}\right| < \frac{1}{2q_n^2}
  \text{ or }
  \left| \alpha - \frac{p_{n+1}}{q_{n+1}}\right| < \frac{1}{2q_{n+1}^2}
  $$
(moreover, for every positive integer $n$, there is $k\in\{n-1,n,n+1\}$ with 
$|\alpha-\frac{p_k}{q_k}|<\frac{1}{\sqrt{5}q_n^2}$).

\noindent $\bullet$  If $\bigl| \alpha-\frac pq\bigr| < \frac{1}{2q^2}$ then $\frac pq$ is a convergent of the continued fraction of $\alpha$.

This formula for $k(\al)$ implies that we have the following alternative
definition of the Lagrange spectrum $L$:

Let $\Sigma=({\nb^*})^{\zb}$ be the set of all bi-infinite sequences of
positive integers. If $\und\theta=(a_n)_{n\in\zb}\in \Sigma$, let $\al_n=[a_n;a_{n+1},a_{n+2},\dots]$ and $\be_n=[0;a_{n-1},a_{n-2},\dots], \forall n \in \zb$. We define
$f(\und\theta)=\al_0+\be_0=[a_0; a_1,a_2,\dots]+[0; a_{-1}, a_{-2},\dots]$.
We have $$L=\{\limsup_{n\to\infty} f(\sigma^n \und\theta), \und\theta\in\Sigma\},$$ 
where $\sigma\colon\Sigma\to\Sigma$ is the shift defined by
$\sigma((a_n)_{n\in\zb})=(a_{n+1})_{n\in\zb}$.

Let us define the Markov spectrum $M$ by 
$$M=\{\sup_{n\in\zb} f(\sigma^n 
\und\theta), \und\theta\in\Sigma\}.$$ 
It also has an arithmetical interpretation, namely
$$
M=\{(\inf_{(x,y)\in{\zb}^2\setminus(0,0)}  |f(x,y)|)^{-1},\quad
f(x,y)=a x^2 + bxy+cy^2, \quad b^2-4ac=1\}.
$$
It is well-known (see \cite{CF}) that $M$ and $L$ are closed sets of the real line
and $L\subset M$.

We have the following result about the Markov and Lagrange spectra:

\bigskip

\begin{theorem}(\cite{M3}; see also \cite{Mat})
Given $t\in\re$ we have
$$
HD(L\cap(-\infty,t))= HD(M\cap(-\infty,t))=: d(t)
$$
and $d(t)$ is a continuous surjective
function from $\re$ to $[0,1]$. Moreover:

i) $d(t)=\min\{1,2 D(t)\}$, where
$D(t):=HD(k^{-1}(-\infty,t)) = HD(k^{-1}(-\infty,t])$ is a continuous function from $\re$ to $[0,1)$.

ii) $\max\{t\in\re\mid d(t)=0\}=3$

iii) $d(\sqrt{12})=1$.
\end{theorem}

A fundamental tool in the proof of this result is the theorem below.

We say that a $C^2$-regular Cantor set on the real line is {\it essentially affine} if there is a $C^2$ change of coordinates for which the dynamics that defines the corresponding Cantor set has zero second derivative on all points of that Cantor set. Typical $C^2$-regular Cantor sets are not essentially affine. 

The {\it scale recurrence lemma}, which is the main technical lemma of \cite{MY}, can be used in order to prove the following

\begin{theorem}(\cite{M4})
If $K$ and $K'$ are regular Cantor sets of class $C^2$ and $K$ is non essentially affine, then $HD(K+K')=\min\{HD(K)+HD(K'), 1\}.$
\end{theorem}

\begin{remark}
There is a presentation of a version of this result (with a slightly different hypothesis) in \cite{S}. That version is also proved by Hochman and Shmerkin in \cite{HS}. 
\end{remark}
 
The results of Markov, Hall and Freiman mentioned above imply that the Lagrange and Markov spectra coincide below $3$ and above $c_F$. Nevertheless, Freiman (\cite{F1}) showed in 1968 that $M\setminus L\neq\emptyset$ by exhibiting a number $\sigma\simeq 3.1181\dots\in M\setminus L$.

In 1973, Freiman \cite{F2} showed that
$$\alpha_{\infty}:=\lambda_0(A_{\infty}):=[2; \overline{1_2, 2_3, 1, 2}] + [0; 1, 2_3, 1_2, 2, 1, \overline{2}]\in M \setminus L$$
In a similar vein, Theorem 4 in Chapter 3 of Cusick-Flahive book \cite{CF} asserts that
$$\alpha_n:=\lambda_0(A_n):= [2; \overline{1_2, 2_3, 1, 2}] + [0; 1, 2_3, 1_2, 2, 1, 2_n, \overline{1, 2, 1_2, 2_3}]\in M\setminus L$$
for all $n\geq 4$. In particular, $\alpha_{\infty}$ is not isolated in $M\setminus L$.

In collaboration with C. Matheus, we proved the following results about $M\setminus L$: 

Let $X$ be the Cantor set
\begin{equation}\label{e.Cantor-X}
X:=\{[0;\gamma]:\gamma\in\{1,2\}^{\mathbb{N}} \textrm{ not containing the subwords in } P\}
\end{equation}
where
$$P:=\{21212, 2121_3, 1_3212, 12121_2, 1_22121, 2_3121_22_21, 12_21_2212_3, 12_3121_22_2,
2_21_2212_31\}$$
Also, let
$$b_{\infty} := [2;\overline{1_2,2_3,1,2}]+[0;\overline{1,2_3,1_2,2}] = \frac{\sqrt{18229}}{41} = 3.2930442439\dots$$
and
\begin{eqnarray*}
B_{\infty} &:=& [2;1,\overline{1,2_3,1,2,1_2,2,1_2,2}]+[0;1,2_3,1_2,2,1,2_3,1_2,2,1,2_2,\overline{1,2_3,1,2,1_2,2,1_2,2}] \\
&=& 3.2930444814\dots
\end{eqnarray*}

\begin{theorem}\label{t.M-L-piece-HD}  (\cite{MM1}) $\{b_{\infty}, B_{\infty}\}\subset L$, $L\cap (b_{\infty}, B_{\infty})=\emptyset$ and  $HD((M\setminus L)\cap (b_{\infty}, B_{\infty})) = HD(X)$.
\end{theorem}

We implemented the algorithm of Jenkinson-Pollicott (see \cite{JP16}) and we obtained the \emph{heuristic} approximation $HD(X)=0.4816\cdots$. We also we exhibited a Cantor set $K(\{1,2_2\})\subset X$ whose Hausdorff dimension can be easily (and rigorously) estimated as $0.353<HD(K(\{1,2_2\}))<0.35792$ via some classical arguments explained in Palis-Takens book \cite{PT1}.

By exploiting the arguments establishing the above Theorem, we are able to exhibit new numbers in $M\setminus L$, including a constant $c\in M\setminus L$ with $c>\alpha_4$:

\begin{proposition}\label{p.new-numbers} The largest element of $(M\setminus L)\cap (b_{\infty}, B_{\infty})$ is
$$c = \frac{77+\sqrt{18229}}{82}+\frac{17633692-\sqrt{151905}}{24923467}=3.29304447990138\dots$$
To our best knowledge, $c$ is the largest known element of $M\setminus L$.
\end{proposition}

We also proved (\cite{MM2}) that $M\setminus L$ \emph{doesn't} have full Hausdorff dimension: 

\begin{theorem}\label{t.prrova}
$HD(M\setminus L) < 0.986927$.
\end{theorem}

One can get better \emph{heuristic} bounds for $HD(M\setminus L)$ thanks to the several methods in the literature to numerically approximate the Hausdorff dimension of Cantor sets of numbers with prescribed restrictions of their continued fraction expansions. By implementing the ``thermodynamical method'' introduced by Jenkinson--Pollicott in \cite{JP01}, we obtained the heuristic bound $HD(M\setminus L)<0.888$. 

Our proof of Theorem \ref{t.prrova} relies on the control of several portions of $M\setminus L$ in terms of the sum-set of a Cantor set associated to continued fraction expansions prescribed by a ``symmetric block'' and a Cantor set of irrational numbers whose continued fraction expansions live in the ``gaps'' of a ``symmetric block''. As it turns out, such a control is possible thanks to our key technical Lemma saying that a sufficiently large Markov value given by the sum of two continued fraction expansions systematically meeting a ``symmetric block'' must belong to the Lagrange spectrum. 

It follows that $M\setminus L$ has empty interior, and so, since $M$ and $L$ are closed subsets of $\mathbb{R}$, $\overline{int(M)}=\overline{int(L)}\subset L\subset M$. In particular, we have the following
\begin{corollary}
$int(M)=int(L)$.
\end{corollary}
As a consequence, we recover the fact, proved in \cite{F3}, that the biggest half-line contained in $M$ coincides with the biggest half-line $[c_F,\infty)$ contained in $L$.

\subsection{Other results on the fractal geometry of Diophantine approximations}

In collaboration with Y. Bugeaud (\cite{BM1}), we proved some results on sets of exact approximation order by rational numbers:

For a function $\Psi : (0,+\infty) \to (0,+\infty)$, let
$$
{\cal K}(\Psi) := \biggl\{ \xi \in \R  : \biggl| \xi - 
{p \over q} \biggr| < {\Psi(q)} \ens
{\hbox{ for infinitely many rational numbers }} {p \over q}
\biggr\}
$$
denote the set of $\Psi$-approximable real numbers and let
$$
{\rm Exact}(\Psi)
:= {\cal K}(\Psi) \setminus \bigcup_{m \ge 2} \, 
{\cal K} \bigl( (1-1/m) \Psi \bigr)
$$
be the set of real numbers approximable 
to order $\Psi$ and to no better 
order.

The {\it lower order at infinity} $\lambda(g)$ of
a function $g : (0,+\infty) \to (0,+\infty)$ is defined by
$$
\lambda(g) = \liminf_{x \to + \infty} \, {\log g(x) \over \log x}.
$$

We say that a function
$\Psi : (0,+\infty) \to (0,+\infty)$ satisfies assumption $(*)$ if
\sm
{\it $\Psi(x) = o(x^{-2})$
and there exist real numbers $c$, $\tilde c$ 
and $n_0$ with $1 \le \tilde c < 4$ such that, 
if the positive integers $m$, $n$ satisfy 
$m > n \ge n_0$, then $\Psi(m) m^c \le \tilde c \Psi(n) n^c$.}
\sm
We emphasize that the real number $c$ occurring in $(*)$
may be negative.

\begin{theorem}
Let $\Psi : (0,+\infty) \to (0,+\infty)$ be a function
satisfying assumption $(*)$.
Then the set ${\rm Exact} (\Psi)$ is uncountable.
\end{theorem}

\begin{theorem}
Let $\Psi : (0,+\infty) \to (0,+\infty)$ be a function
satisfying assumption $(*)$. If $\lambda$
denotes the lower order at infinity 
of the function $1/\Psi$, then
$$
\dim {\rm Exact} (\Psi) = \dim {\cal K}(\Psi) = {2 \over \lambda}.
$$
\end{theorem}

In another work collaboration with Y. Bugeaud (\cite{BM2}), we proved that there are no typical real numbers from the point of view of Diophantine approximations, in a sense that we describe in what follows. Let $\Psi$ be an application from the set of positive integers into the set of nonnegative real numbers. Khintchine established that, 
if the function $q \mapsto q^2 \Psi(q)$ is non-increasing 
and the series $\sum_{q \ge 1} q \Psi(q)$ diverges, then the
set $\cK(\Psi)$ has full Lebesgue measure (Beresnevich, Dickinson and Velani proved later the same result assuming that $\Psi$ is just non-increasing). We show that, for almost every real number $\alpha$, there is a function $\Psi$ which satisfies good ``regularity" conditions (on the speed of decreasing of $\Psi$) - for instance $q \mapsto q^2 \Psi(q)$ is non-increasing, such that the series $\sum_{q \ge 1} q \Psi(q)$ diverges but the inequality $|\alpha - p/q| < \Psi(q)$ has no rational solution $p/q$. 

Khintchine also showed that if the series $\sum_{q \ge 1} q \Psi(q)$ converges, then the set $\cK(\Psi)$ has zero Lebesgue measure. We show that, for almost every real number $\alpha$, there is a function $\Psi$ which satisfies good ``regularity" conditions (for instance $q \mapsto q^2 \Psi(q)$ is non-increasing), such that the series $\sum_{q \ge 1} q \Psi(q)$ converges but the inequality $|\alpha - p/q| < \Psi(q)$ has infinitely many rational
solutions $p/q$.

We also compute Hausdorff dimensions of sets of exceptions to our results (in terms of the regularity conditions on $\Psi$):

\begin{theorem}
Let $\lambda$ be a real number in the interval $[1,2]$. 
We define $X_{\lambda}$ the set of irrational numbers $\xi$ such that, for every function 
$\Psi$ with $q \mapsto q^{\lambda}\Psi(q)$ non-increasing and satisfying 
$$
\sum_{q=1}^{\infty} \, q \Psi(q) = + \infty,  
$$
the inequality
$$
\biggl| \xi - {p \over q} \biggr| < {\Psi(q)}   
$$
has infinitely many rational solutions.
We proved that $X_1$ is empty and that, for $\lambda\in (1,2]$, 
the Hausdorff dimension of $X_{\lambda}$
is $\lambda/2$ (and so, when $\lambda$ varies, assume all values in the interval $(1/2,1]$).
\end{theorem}

\begin{theorem}
Let $b$ be a positive real number. We define $Y_b$ as the set
of the irrational numbers $\xi$ such that, for every function 
$\Psi$ with$q \mapsto q (\log \log q)^{b \log \log q} \Psi(q)$ 
non-increasing and satisfying
$$
\sum_{q=1}^{\infty} \, q \Psi(q) = + \infty,  
$$
the inequality
$$
\biggl| \xi - {p \over q} \biggr| < {\Psi(q)}   
$$
has infinitely many rational solutions.
Then the Hausdorff dimension of $Y_b$
is ${1 \over {1+ {\rm e}^{1/b}}}$ (and so, when $b$ varies, assume all values in the interval $(0,1/2)$).
\end{theorem}

One of the tools we used in the proof of this last theorem is the results from \cite{L}: given $b, c>1$ define $\tilde\Xi(b,c)=\{\xi=[0;a_1,a_2,...] \in [0,1] : a_n\ge c^{b^n} \ 
\hbox{for every $n \ge 1$}\}$ and $\Xi(b,c)=\{\xi=[0;a_1,a_2,...] \in [0,1] : a_n\ge c^{b^n} \ 
\hbox{for infinitely many $n \ge 1$}\}$. Then $HD(\tilde\Xi(b,c))=HD(\Xi(b,c))={1 \over 1+b}$. In one direction, a more precise result is proved in \cite{FWLT}, according to which the Hausdorff dimension of the set
$$
K:=\{\xi=[0;a_1,a_2,...] \in [0,1] : \exp({\tilde b}^n)\le a_n \le 3 \exp({\tilde b}^n), \ 
\hbox{for every $n \ge 1$}\}
$$
is equal to
${1 \over 1+\tilde b}$. We notice that an older result from \cite{G} states that the Hausdorff dimension of the set $\{\xi=[0;a_1,a_2,...] \in [0,1] : \lim a_n=\infty\}$ is $1/2$.

\section{Back to Dynamical Systems (and Differential Geometry)}

As we have seen, the sets $M$ and $L$ can be defined in terms of symbolic dynamics. Inspired by these characterizations, we may associate to a dynamical system together with a real function generalizations of the Markov and Lagrange spectra as follows:
\begin{definition}	
	Given a map $\psi:X \to X$ and a function $f:X\to \mathbb R$, we define the associated dynamical Markov and Lagrange spectra as \hfill \break 
	$M(f,\psi)=\{\text{sup}_{n\in \nb}f(\psi^n(x)), x\in X\}$ and \hfill \break
	$L(f,\psi)=\{\text{limsup}_{n\to\infty}f(\psi^n(x)), x\in X\}$, respectively.
	
	Given a flow $(\varphi^t)_{t\in \mathbb R}$ in a manifold $X$, we define the associated dynamical Markov and Lagrange spectra as 
	$M(f,(\varphi^t))=\{\text{sup}_{t\in \R}f(\varphi^t(x)), x\in X\}$ and \hfill \break
	$L(f,(\varphi^t))=\{\text{limsup}_{t\to\infty}f(\varphi^t(x)), x\in X\}$, respectively.
\end{definition}
\bigskip
In a work in collaboration with A. Cerqueira and C. Matheus (\cite{CMM}), we prove the following result, which generalizes a corresponding fact in the context of the classical Markov and Lagrange spectra:
\begin{lemma}
Let $(\varphi, f)$ be a generic pair, where $\varphi \colon M^2\to M^2$ is a diffeomorphism with $\Lambda\subset M^2$ a hyperbolic set for $\varphi$ and $f:M\to \mathbb R$ is $C^2$. Let $\pi_s, \pi_u$ be the projections of the horseshoe $\Lambda$ to the stable and unstable regular Cantor sets $K^s, K^u$ associated to it (along the unstable and stable foliations of $\Lambda$). Given $t\in \mathbb R$, we define 
$$\Lambda_t=\bigcap_{m\in \mathbb Z}\varphi^m(\{p\in \Lambda|f(p)\le t\}),$$
$$K_t^s=\pi_s(\Lambda_t), K_t^u=\pi_u(\Lambda_t).$$
Then the functions $d_s(t)=HD(K_t^s)$ and $d_u(t)=HD(K_t^u)$ are continuous and coincide with the corresponding box dimensions.
\end{lemma}

The following result is a consequence of the scale recurrence lemma of \cite{MY} (see \cite{M4}):
\begin{lemma}
Let $(\varphi, f)$ be a generic pair, where $\varphi \colon M^2\to M^2$ is a diffeomorphism with $\Lambda\subset M^2$ a hyperbolic set for $\varphi$ and $f:M\to \mathbb R$ is $C^2$. Then
$$HD(f(\Lambda))=\min(HD(\Lambda),1).$$
Moreover, if $HD(\Lambda)>1$ then $f(\Lambda)$ has persistently non-empty interior.
\end{lemma}

Using the previous lemmas we prove a generalization of the results on dimensions of the dynamical spectra:
\begin{theorem}(\cite{CMM})
Let $(\varphi, f)$ be a generic pair, where $\varphi \colon M^2\to M^2$ is a {\bf conservative} diffeomorphism with $\Lambda\subset M^2$ a hyperbolic set for $\varphi$ and $f:M\to \mathbb R$ is $C^2$. Then
$$HD(L(f,\Lambda)\cap(-\infty,t))=HD(M(f,\Lambda)\cap(-\infty,t))=:d(t)$$
is a continuous real function whose image is $[0,\min(HD(\Lambda),1)]$.
\end{theorem}\
(here we use the notation $L(f,\Lambda):=L(f,\varphi|_{\Lambda})$).

In \cite{HMU}, P. Hubert, L. Marchese and C. Ulcigrai introduced in a similar way to the above generalizations the Lagrange spectra of closed-invariant loci for the action of $SL(2,\re)$ on the moduli space of translation surfaces, in the context of Teichm\"uller dynamics, and proved that several of these spectra contain a Hall's ray.

Moreira and Roma\~na prove the following result on Markov and Lagrange spectra for horseshoes:

\begin{theorem}[\cite{MR1}]
Let $\Lambda$ be a horseshoe associated to a $C^2$-diffeomorphism $\varphi$ such that  $HD(\Lambda)>1$. Then there is, arbitrarily close to $\varphi$ a diffeomorphism $\varphi_{0}$ and a $C^{2}$-neighborhood $W$ of $\varphi_{0}$ such that, if $\Lambda_{\psi}$ denotes the continuation of $\Lambda$ associated to $\psi\in W$, there is an open and dense set $H_{\psi}\subset C^{1}(M,\re)$ such that for all $f\in H_{\psi}$, we have 
\begin{equation*}
int \ L(f,\Lambda_{\psi})\neq\emptyset \ \text{and} \ int  \ M(f,\Lambda_{\psi})\neq\emptyset,
\end{equation*}
where $int \, A$ denotes the interior of $A$.
\end{theorem}

Recently, D. Lima proved that, for typical pairs $(f,\Lambda)$ as in the above theorem, 
$$\sup\{t\in\re|HD(M(f,\Lambda)\cap(-\infty,t)<1\}=\inf\{t\in\re|int(L(f,\Lambda)\cap(-\infty,t)\ne\emptyset\},$$
and Moreira (\cite{M5}) proved that, for typical pairs $(f,\Lambda)$ as above, the minima of the corresponding Lagrange and Markov dynamical spectra coincide and are given by the image of a periodic point of the dynamics by the real function, solving a question by Yoccoz.

The classical Markov and Lagrange spectra can also be characterized as sets of maximum heights and asymptotic maximum heights, respectively, of geodesics in the modular surface $N=\mathbb{H}^{2}/PSL(2,\mathbb{Z})$. Moreira and Roma\~na extend in \cite{MR2} the fact that these spectra have non-empty interior to the context of negative, non necessarily constant curvature as follows:

\begin{theorem}(\cite{MR2})
Let $M$ provided with a metric $g_0$ be a complete noncompact surface $M$ with finite Gaussian volume and Gaussian curvature bounded between two negative constants, \emph{i.e.}, if $K_{M}$ denotes the Gaussian curvature, then there are constants $a,b > 0$ such that
$$-a^2\leq K_{M}\leq-b^2<0.$$
Denote by $SM$ its unitary tangent bundle and by $\phi$ its geodesic flow.

Then there is a metric $g$ close to $g_0$ and a dense and $C^{2}$-open subset $\mathcal{H}\subset C^{2}(SM,\re)$  such that 
$$int\ M(f,\phi_{g})\neq \emptyset \ \text{and} \ int\ L(f,\phi_{g})\neq \emptyset$$
for any $f\in \mathcal{H}$,
where $\phi_{g}$ is the vector field defining the geodesic flow of the metric $g$. 

Moreover, if $X$ is a vector field sufficiently close to $\phi_g$ then 
$$int\ M(f,X)\neq \emptyset \ \text{and} \ int\ L(f,X)\neq \emptyset$$
for any $f\in \mathcal{H}$.
\end{theorem}

Recently, in collaboration with Pac\'ifico and Roma\~na (\cite{MPR}), we proved an analogous result for Lorenz flows. Roma\~na also proved (\cite{R}) a corresponding result for Anosov flows in compact $3$-dimensional manifolds.

Combining the techniques of this result with those of the above results in collaboration with A. Cerqueira and C. Matheus, Cerqueira, Moreira and Roma\~na proved the following

\begin{theorem}(\cite{CMR})
Let $(\varphi, f)$ be a generic pair, where $\varphi \colon M^2\to M^2$ is a {\bf conservative} diffeomorphism with $\Lambda\subset M^2$ a hyperbolic set for $\varphi$ and $f:M\to \mathbb R$ is $C^2$. Then
$$HD(L(f,\Lambda)\cap(-\infty,t))=HD(M(f,\Lambda)\cap(-\infty,t))=:d(t)$$
is a continuous real function whose image is $[0,\min(HD(\Lambda),1)]$.
\end{theorem}\

%\frenchspacing

\end{document}